\documentclass[10pt]{icsnews}

\newslettername{ICS News}
\issuedate{May 2026}

\graphicspath{{images/}}

\makeatletter
\renewcommand{\ResearchHighlightSection}[2]{%
  \ResetArticleCounters
  \hypertarget{researchhighlight#1}{}%
  \begin{center}
    {\Large\textbf{#2\footnotemark}}%
  \end{center}%
  \footnotetext{This expository article will appear as an invited Research Highlight in the 2026 INFORMS Computing Society Newsletter.}%
  \vspace{2pt}%
}
\renewcommand{\NewsAuthor}[1]{%
  \begin{center}
  #1\vspace{10pt}
  \par\normalfont
  \end{center}%
}
\makeatother

\begin{document}


\ResearchHighlightSection{1}{Stepsize Hedging: \\ An Alternative Mechanism for Accelerating Gradient Descent}
\NewsAuthor{%
\begin{tabular}{cc}
Jason M. Altschuler & Pablo A. Parrilo \\
UPenn & LIDS - MIT
\end{tabular}%
}




\newcommand*{\eps}{\varepsilon}
\newcommand*{\R}{\mathbb{R}}
\newcommand{\note}[1]{{\textcolor{red}{[#1]}}}
\newcommand{\todo}{\note{TODO}}
\newcommand{\todoc}{\note{[citation]}}
\newcommand{\ja}[1]{\note{JASON: #1}}
\newcommand{\pp}[1]{\note{PABLO: #1}}
\providecommand{\BIBand}{and}


\section{Introduction}\label{sec:intro}

Consider the following question that arises in every introductory course on convex optimization:
\begin{align*}
    \text{\emph{What is the optimal choice of stepsizes for gradient descent?}} 
\end{align*}
Gradient descent (GD) is arguably the most foundational algorithm for continuous optimization and has a commensurately extensive literature. GD was first proposed by Cauchy in the 1800s and remains the workhorse in large-scale optimization applications across engineering, data science, and artificial intelligence. Yet, the fundamental question stated above---how to choose the sole parameter in GD---remained open even in seemingly simple convex settings, highlighting the potential for untapped algorithmic opportunities in those foundational settings and beyond.

\textbf{Short or long stepsizes?} Let us begin by recalling the definition of GD:
\begin{align}
    x_{t+1} = x_t - \alpha_t \nabla f(x_t)\,.
\end{align}
The basic idea is that the negative gradient is a descent direction, meaning that if one moves a small amount in this direction, then the value of the objective $f$ decreases. This follows immediately by a Taylor expansion. The algorithmic question at hand is: \textit{how far} should one move? A small stepsize $\alpha_t$ guarantees progress, but the progress will be correspondingly small. A large stepsize $\alpha_t$ is a natural idea to get more mileage out of the descent direction, but this may overshoot the optimum. 

\textbf{Mainstream prescription for smooth convex optimization: constant stepsize schedules.} A natural way to quantify this tradeoff is to optimize the progress made in a single iteration. In the standard setup of smooth convex optimization, this one-step question has a well-known answer:
\begin{gather}
	\text{\emph{For one iteration of GD, there is a stepsize} } \bar{\alpha} \text{ \emph{that achieves the fastest worst-case convergence rate}}.
	\label{eq:intro:1-step}
\end{gather}
This stepsize is explicit.\footnote{For example, $\bar{\alpha} = 1/M$ for $M$-smooth convex objectives, and $\bar{\alpha} = 2/(M+m)$ for objectives that are also $m$-strongly convex.} This one-step guarantee underlies the textbook prescription: use the constant schedule $\alpha_t = \bar{\alpha}$ at every iteration $t$. See for example the textbooks~\cite{bubeck-book,boyd,polyakbook,hazan2016introduction,nesterov-survey,luenberger1984linear,BertsekasNonlinear}. 

However, even after optimizing $\bar{\alpha}$, this constant stepsize schedule may lead to slow convergence. This has motivated an extensive literature on \emph{accelerated} first-order methods which famously modify GD by adding momentum or other internal dynamics, see the survey~\cite{d2021acceleration}. Many alternative stepsize schedules have also been proposed---for example exact line search, Armijo-Goldstein rules, Polyak-type schedules, and Barzilai-Borwein-type schedules---but none had led to an analysis that improves over the textbook constant-stepsize rate. The resulting conventional wisdom was that GD cannot be accelerated merely by changing its stepsizes.

\textbf{Faster convergence via hedged stepsizes?} Is this conventional wisdom correct? Surprisingly, the answer is no. The key observation is that optimality of the stepsize $\bar{\alpha}$ for one iteration does \emph{not} imply optimality of repeating the stepsize $\bar{\alpha}$ for multiple iterations. The stepsize $\bar{\alpha}$ is chosen to protect against the worst case for a single iteration. But when GD is run for multiple iterations, the bad instances for different stepsizes need not be compatible with one another. A short step may be too conservative on one kind of instance (flat objectives $f$), while a long step may overshoot on another (steep objectives $f$); over several iterations, these opposing weaknesses can fail to align (since convex objectives $f$ are rigid and cannot change curvature arbitrarily). This suggests a multi-step opportunity:
\begin{gather}
	\text{\emph{Can one combine stepsizes that are individually suboptimal}} \nonumber \\
	\text{\emph{to obtain faster convergence over many iterations?}}
	\label{eq:intro:hedge}
\end{gather}
We refer to this idea as \emph{stepsize hedging}: rather than optimizing each step in isolation, choose a time-varying schedule that hedges between different worst-case behaviors.

\paragraph*{Motivation: the special case of quadratics.} Time-varying stepsizes were initially explored (only) in the special case of convex quadratic optimization. A classic result due to Young in 1953~\cite{young53} shows that in this setting, the optimal stepsizes are non-constant and related to the roots of Chebyshev polynomials. This elegant result is recalled in detail in \S\ref{sec:quadratic}. However, many core phenomena in the quadratic setting do not extend beyond. For example, these Chebyshev stepsizes can make GD diverge beyond the special case of quadratics. Over the past 70 years, the continuous optimization community has devoted significant effort to developing alternative stepsize schedules beyond quadratic optimization. These are often empirically helpful. However, despite significant effort and many candidate stepsize schedules (even adaptive), there was no theoretical evidence that \emph{any} stepsize schedule could lead to \emph{any} speedup over the textbook GD convergence rate for (non-quadratic) convex optimization.

\subsection{Main result}\label{ssec:intro:result}

In the past decade, a line of work has shown that such an algorithmic opportunity persists: time-varying stepsizes can accelerate GD for (non-quadratic) convex optimization. The first such result was shown in the thesis~\cite{altschuler2018greed} of the first author (advised by the second author), and an exciting flurry of ensuing work has sought to push this algorithmic opportunity to its limit. See the discussion of related work in \S\ref{ssec:intro:related}. 

\par We summarize here the results of the two papers~\cite{hedging1,hedging2}. In these papers, we propose an unconventional stepsize schedule and show that it accelerates GD for smooth convex optimization. We term this the 
\emph{silver stepsize schedule} due to the occurrence of the \emph{silver ratio} $\rho = 1 + \sqrt{2}$. It is explicit and non-adaptive, see \S\ref{sec:silver-stepsizes} for the definition and a full discussion. Here we highlight the main results, namely the accelerated convergence rates that this enables for the convex and strongly-convex settings, respectively.

\par Below, for shorthand, we say that an iterative algorithm initialized at $x_0$ minimizes an $M$-smooth convex function $f$ to $\eps$ error if its final iterate $x_n$ satisfies $\tfrac{f(x_n) - f(x^*)}{M\|x_0 - x^*\|^2} \leq \eps$, where $x^*$ denotes a minimizer of $f$. The factor of $M\|x_0 - x^*\|^2$ normalizes the performance to make it scale-invariant and meaningful.

\begin{theorem}[Silver stepsizes for convex optimization~\cite{hedging2}]\label{thm:convex}
    Fix any accuracy $\eps$ and smoothness parameter $M$. There exists a stepsize schedule of length
    \begin{align*}
        n \asymp \eps^{- \log_{\rho} 2} \approx \eps^{-0.7864}  \end{align*}
    such that GD $\eps$-minimizes any $M$-smooth convex function $f : \R^d \to \R$ in any dimension $d$.
\end{theorem}

\begin{theorem}[Silver stepsizes for strongly convex optimization~\cite{hedging1}]\label{thm:strongly-convex}
    Fix any accuracy $\eps$, strong convexity parameter $m$, and smoothness parameter $M$. There exists a stepsize schedule of length
    \begin{align*}
        n \asymp \kappa^{\log_{\rho} 2} \log \frac{1}{\eps} \approx \kappa^{0.7864} \log \frac{1}{\eps}
    \end{align*}
     such that GD $\eps$-minimizes any $m$-strongly convex, $M$-smooth function $f : \R^d \to \R$ in any dimension $d$. Here $\kappa = M/m \geq 1$ denotes the condition number.
\end{theorem}

For comparison, Theorem~\ref{thm:convex} improves the textbook rate of constant-step GD from $O(\eps^{-1})$ to roughly $O(\eps^{-0.7864})$. In the strongly convex setting, Theorem~\ref{thm:strongly-convex} analogously improves the textbook rate of GD from $O(\kappa \log \tfrac{1}{\eps})$ to roughly $O(\kappa^{0.7864} \log \tfrac{1}{\eps})$. These results establish stepsize hedging as an alternative mechanism for accelerating GD. Interestingly this does not match the fully accelerated rate obtained by modifying GD beyond stepsizes~\cite{nesterov-agd}. Our asymptotic rates are conjectured optimal among all (non-adaptive\footnote{A non-adaptive schedule is a predetermined sequence $\alpha_0,\alpha_1,\ldots$ that may depend on problem parameters such as $M,m,n$, but not on the observed iterates or gradients. In contrast, an adaptive schedule may choose $\alpha_t$ using information from earlier iterations.}) stepsize schedules.

\subsection{Related work}\label{ssec:intro:related}

There is an extensive literature on accelerated first-order algorithms, e.g.,~\cite{d2021acceleration,nesterov-survey,nem-yudin}. As opposed to the mainstream approach which modifies GD, we focus here on an alternative mechanism for accelerating GD: time-varying stepsize schedules.

As mentioned above, the benefit of time-varying stepsizes was classically known for the special case of quadratic optimization since Young in 1953~\cite{young53}, but remained open beyond quadratics. The first such results were by~\citet{altschuler2018greed}, who showed that time-varying hedging can lead to improved rates by analyzing $n=2,3$ in the strongly convex setting.~\citet{daccache2019performance} and~\citet{eloi2022worst} exhaustively extended these $n=2,3$ results to related settings and performance metrics. An exciting line of subsequent work numerically optimized longer schedules. Hedging over larger horizons $n$ further improves the convergence rate, but searching for optimal stepsizes is a non-convex problem that is computationally difficult as $n$ increases.~\citet{gupta22} developed a branch-and-bound framework to compute good schedules up to $n=50$ for smooth convex optimization.~\citet{grimmer23} developed a technique to round these branch-and-bound solutions to exact rational certificates; this allowed him to extend these approximate stepsize schedules up to $n=127$ in order to get a larger constant-factor improvement. 

\par How far can this be pushed? As $n$ increases (beyond constants), can this time-varying hedging lead to \emph{asymptotic} acceleration? What are the optimal stepsizes? Intriguing fractal-like patterns were observed for small $n$; do these qualitative behaviors persist for large $n$?

The two papers~\cite{hedging1,hedging2} highlighted in this note answered this asymptotic-acceleration question affirmatively by introducing the \emph{silver stepsize schedule}. This schedule is built recursively from the same short-step/long-step splitting pattern that appears in the optimal two-step schedule of~\cite{altschuler2018greed}. In
the smooth convex setting (Theorem~\ref{thm:convex}), the silver stepsizes improve the iteration complexity of GD from the textbook rate $O(\eps^{-1})$ to $O(\eps^{-\log_{\rho} 2}) \approx O(\eps^{-0.78})$. In the $m$-strongly convex and $M$-smooth setting (Theorem~\ref{thm:strongly-convex}), they similarly improve the dependence on the condition number $\kappa=\tfrac{M}{m}$ from $O(\kappa \log \tfrac{1}{\eps})$ for constant stepsizes to $O(\kappa^{\log_{\rho} 2}\log \tfrac{1}{\eps})
    \approx O(\kappa^{0.78}\log \tfrac{1}{\eps})$. 
These asymptotic rates were conjectured to be asymptotically optimal among all non-adaptive stepsize schedules~\cite{hedging1,hedging2}. (Concurrent work by~\citet{GrimmerShuWang} proved asymptotic acceleration using a different stepsize schedule, but at a suboptimal convergence rate, with exponent $\approx 0.95$ rather than $\log_{\rho}2\approx 0.78$.)

Excitingly, in only the few years since, many papers have followed up; see \S\ref{sec:outlook} for a discussion of this burgeoning area. 

\begin{table}[t]
\small
	\begin{tabular}{|c|c|c|}
		\hline
		& \textbf{Quadratic }                          & \textbf{Convex      }                              \\ \hline
		\textbf{Mainstream stepsizes}  & $\Theta(\kappa)$ by constant stepsizes (folklore)  & $\Theta(\kappa)$ by constant stepsizes (folklore)   \\
        \textbf{Optimized stepsizes}    & $\Theta(\sqrt{\kappa})$ by Chebyshev Stepsizes~\cite{young53} & \color{blue}{ $O(\kappa^{\log_{\rho} 2})$ by Silver Stepsizes~\cite{hedging1} (Theorem~\ref{thm:strongly-convex})} \\ 
		\textbf{Additional dynamics} & $\Theta(\sqrt{\kappa})$ by Heavy Ball~\cite{polyak1964some}          & $\Theta(\sqrt{\kappa})$ by Nesterov Acceleration~\cite{nesterov-agd}     \\ \hline
	\end{tabular}
	\caption{
    Iteration complexity of various approaches for minimizing a $\kappa$-conditioned convex function. The dependence on the accuracy $\eps$ is omitted as it is always $\log 1/\eps$. The story is analogous without strong convexity: rates of the form $O(\kappa^a \log 1/\eps)$ here correspond to rates of the form $O(\eps^{-a})$ there. 
    This note focuses on a foundational question: how much mileage can one obtain from optimizing the stepsizes of GD?}
\end{table}

\section{Optimal stepsizes for quadratic optimization (Young 1953)}\label{sec:quadratic}

As mentioned above, in the special case of \emph{quadratic} optimization, it is classically known that time-varying stepsizes can accelerate GD. This result is due to~\citet{young53} and is nowadays taught in introductory optimization courses. We briefly recall this elegant argument below as it provides perhaps the simplest example of the broader algorithmic opportunity. 

\textbf{Young's stepsizes.} Consider running GD on the class $\mathcal{F}$ of quadratic functions $f$ that are $m$-strongly convex and $M$-smooth. What stepsize schedules make GD converge
\begin{align}
    \|x_n - x^*\| \leq R_n\|x_0 - x^*\|
\end{align}
at the fastest possible rate $R_n$?
Without loss of generality after translating, $f(x) = \frac{1}{2}x^T Hx $ where $mI \preceq H \preceq MI$. Since $f$ is quadratic, its gradient is linear $\nabla f(x) = Hx$, hence GD is a linear map $x_{t+1} = x_t - \alpha_t \nabla f(x_t) = (I - \alpha_t H)x_t$, and therefore the $n$-th iterate is
\begin{align}
	x_n = p_n(H)x_0\,,\qquad \text{ where } \qquad p_n(H) = \prod_{t < n}(I - \alpha_t H).
    \label{eq:quad-xn}
\end{align}
Observe that as one ranges over all possible choices of the stepsize schedule $\{\alpha_t\}_{ t < n}$, the polynomial $p_n$ ranges over the set $\mathcal{P}_n$ of all degree-$n$ polynomials satisfying the normalizing condition $p_n(0) = 1$. 
Therefore finding an optimal stepsize schedule is equivalent to finding an optimal polynomial $p_n \in \mathcal{P}_n$. 

\par
\begin{figure}[h]
\centering
\includegraphics[width=0.38\linewidth]{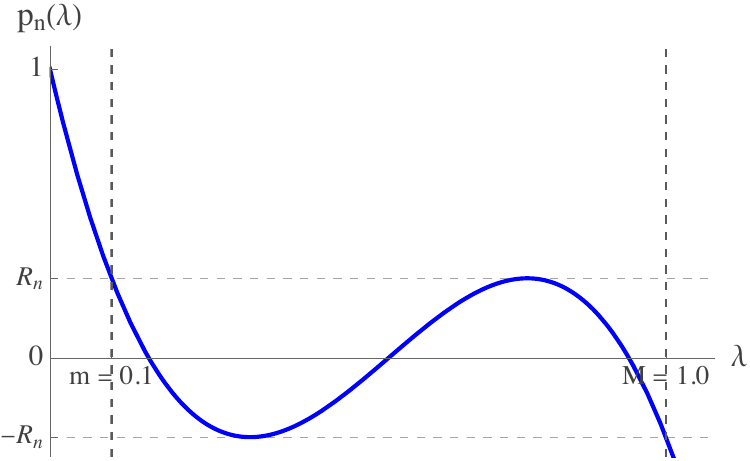}
\caption{ 
   The translated and scaled Chebyshev polynomial $p_n \in \mathcal{P}_n$ has minimal $L_{\infty}$ norm on $[m,M]$. Here visualized for $n=4$, $m=0.1$, $M=1$. In quadratic optimization, the error after GD with stepsizes $\{\alpha_t\}_{t<n}$ is governed by the polynomial $p_n(\lambda)=\prod_{t<n}(1-\alpha_t\lambda)$. Chebyshev stepsizes hedge across curvatures $\lambda \in [m,M]$ by making the worst-case errors $p_n(\lambda)$ equioscillate between $\pm R_n$.
}
\label{fig:chebyshev}
\end{figure}

 What is the optimal polynomial? By the above display and basic properties of the spectral norm,
\begin{align*}
    R_n
    = \sup_{f \in \mathcal{F},\; x_0 \neq x^*} \frac{\|x_n - x^*\|}{\|x_0 - x^*\|}
    =
	\sup_{mI \preceq H \preceq MI, \, x_0 \neq 0} \frac{\|p_n(H)x_0\|}{\|x_0\|}
	=
	\sup_{mI \preceq H \preceq MI}
	\|p_n(H)\|
	=
	\sup_{m \leq \lambda \leq M} 
	|p_n(\lambda)|\,.
\end{align*}
Thus the optimal polynomial $p_n \in \mathcal{P}_n$ is the one with minimal $L_{\infty}$ norm over the interval $[m,M]$. It is classically known that this is the (translated and scaled) Chebyshev polynomial of the first kind, see e.g.,~\cite{Rivlin}. By the definition of $p_n$ in~\eqref{eq:quad-xn}, it follows that the optimal stepsizes $\{\alpha_t\}_{t < n}$  are the inverses of the roots
of this Chebyshev polynomial, in any order. See Figure~\ref{fig:chebyshev} for an illustration and an interpretation of this phenomenon through the lens of \emph{hedging}.

\textbf{Beyond quadratics?} 
This argument establishes that time-varying stepsizes can accelerate GD for quadratic $f$, but does this algorithmic opportunity extend beyond? From an analysis perspective, the argument fails from the very beginning: if $f$ is not quadratic, then GD is not a linear map, and the equivalence to polynomials fails. 
In fact, this issue is not merely an artifact of analysis techniques: fundamental phenomena for the quadratic setting are false beyond. For example, in (non-quadratic) convex optimization, these Chebyshev stepsizes lead to divergence, the stepsize order matters, stepsizes and momentum are not equivalent, etc. As a result, it was unknown if \emph{any} stepsize schedule could lead to \emph{any} improvement over the constant stepsize schedule in the general setting of (non-quadratic) convex optimization. Is this a missed opportunity?

\section{Optimal stepsizes for convex optimization, $n=2$ (Altschuler 2018)}\label{sec:2-step}

The answer is yes. This was first shown in \cite{altschuler2018greed}. In fact, that thesis showed several such results; we state here the simplest one.

Consider $n=2$ steps, the minimal setting where time-varying stepsizes could possibly be advantageous. What two stepsizes $\alpha_0,\alpha_1$ make GD converge
\begin{align}
    \|x_2 - x^*\| \leq R_2 \|x_0 - x^*\|
    \label{eq:2-step-progress}
\end{align}
at the fastest possible rate $R_2$? In the worst case over objectives $f$ that are $m$-strongly convex and $M$-smooth? Obviously $R_2 \leq R_1^2$ is possible by using $\alpha_0 = \alpha_1 = \bar{\alpha}$. Can strictly faster convergence $R_2 < R_1^2$ be achieved? With time-varying stepsizes $\alpha_0 \neq \alpha_1$? Remarkably, the answer is yes:

\begin{theorem}[Theorem 8.10 of \cite{altschuler2018greed}]\label{thm:2-step}
    The stepsizes $\alpha_0, \alpha_1$ that make $R_2$ as small as possible in~\eqref{eq:2-step-progress} are unique,  time-varying ($\alpha_0 \neq \alpha_1$), and strictly improve over the constant stepsize schedule ($R_2 < R_1^2$). 
\end{theorem}

\begin{figure}[h]
\centering
\includegraphics[width=0.47\linewidth]{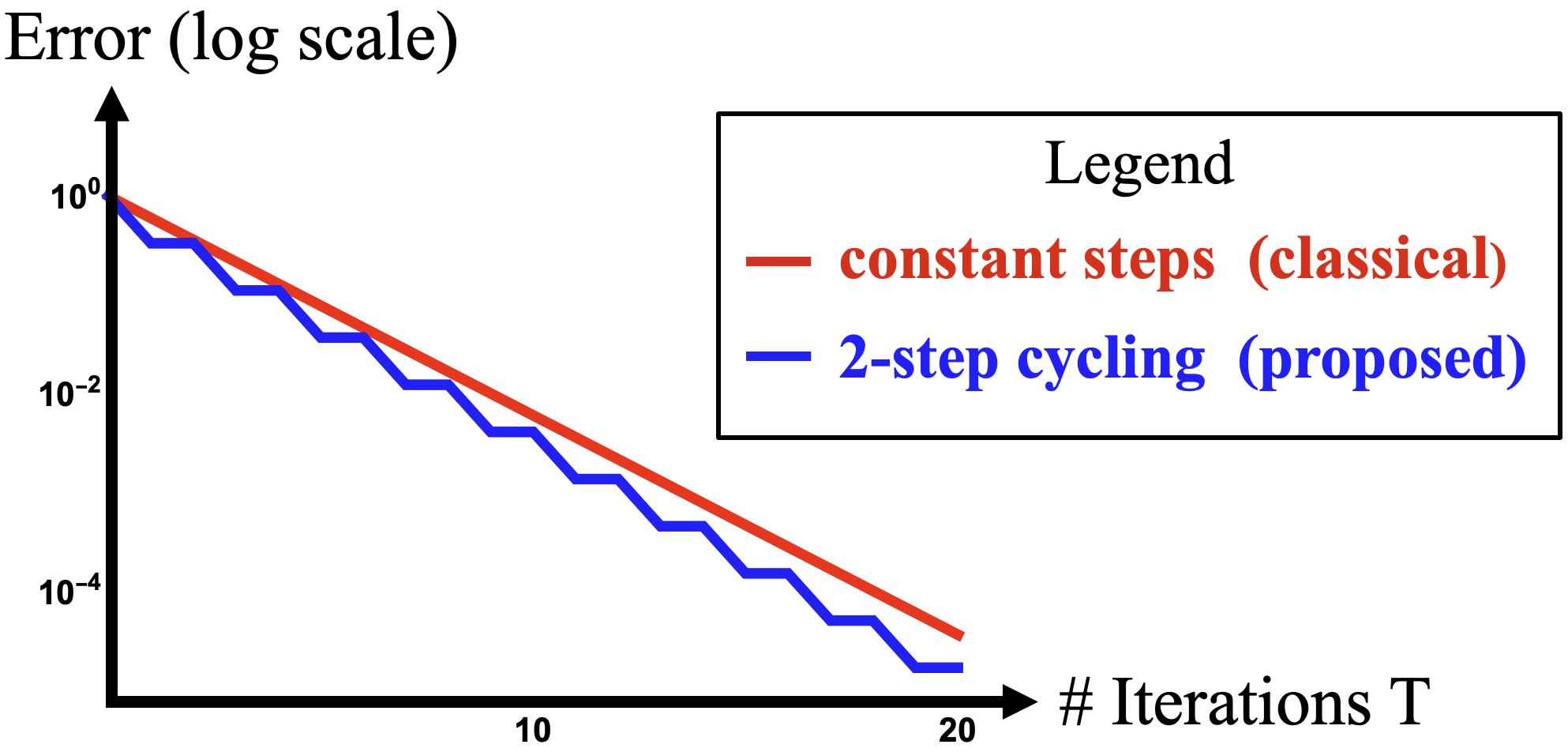}
\caption{ 
    Illustration of Theorem~\ref{thm:2-step}: time-varying stepsizes that cycle between $2$ stepsizes provably speed up GD for (non-quadratic) convex optimization.
    Here consider minimizing functions that are $1$-strongly convex and $4$-smooth. \textbf{Red:} the standard prescription of constant stepsizes $\alpha_t = 2/5$ in each iteration $t$ leads to a convergence rate of $0.6^T$ after $T$ iterations. \textbf{Blue:} better is \emph{cycling} between $\alpha_t = 1/3$ and $1/2$ in even and odd steps; this leads to a faster convergence rate of $
    \smash[t]{(1/\sqrt{3})^T}
    \approx 0.57^T$. 
}
\label{fig:2step-rate}
\end{figure}

As a trivial corollary, cyclically repeating $(\alpha_0, \alpha_1, \alpha_0, \alpha_1, \alpha_0, \alpha_1, \dots)$ strictly improves over the textbook rate for constant stepsizes. Indeed, the two-step convergence~\eqref{eq:2-step-progress} implies $\|x_{2n} - x^*\| \leq R_{2}^n \|x_0 - x^*\|$, which strictly improves over the standard rate $\|x_{2n} - x^*\| \leq R_1^{2n} \|x_0 - x^*\|$ for constant stepsizes. See Figure~\ref{fig:2step-rate}.

The precise values for $\alpha_0, \alpha_1, R_2$ are complicated and not essential for this brief exposition\footnote{For the interested reader, the optimal stepsizes are $\alpha_0 = \frac{2}{m+S}$ and $\alpha_1 = \frac{2}{2M+m-S}$, with corresponding rate $R_2 = \frac{S-M}{2m+S-M}$, where $S = \sqrt{M^2+(M-m)^2}$.}. We only highlight two features of these explicit values that will be built upon later:
\begin{itemize}
    \item \textbf{Stepsize splitting.} The optimal stepsizes $\alpha_0 < \bar{\alpha} < \alpha_1$ are obtained by ``splitting'' the constant stepsize $\bar{\alpha}$ into a shorter step $\alpha_0$ and a longer step $\alpha_1$ as the two roots of a certain explicit quadratic equation.
    \item \textbf{Silver ratio.} The improvement $R_2 < R_1^2$ can be quantified as $R_2 \approx 1 - \tfrac{2(1+\sqrt{2})}{\kappa}$ whereas $R_1^2 \approx 1 - \frac{4}{\kappa}$, in the relevant asymptotic regime as $\kappa \to \infty$. This quantity $\rho = 1 + \sqrt{2}$ is called the silver ratio.
\end{itemize}

Summarizing, Theorem~\ref{thm:2-step} shows a missed algorithmic opportunity: time-varying stepsizes provably make GD converge faster. However, this result only shows a constant factor improvement in the final iteration complexity (from $R_2 < R_1^2$) since it only considers $n=2$ steps. The improvement increases in $n$. But by how much? As $n \to \infty$, what are the optimal stepsizes and rate? How much faster can one accelerate GD?

\section{Silver stepsizes for convex optimization, $n \to \infty$ (Altschuler-Parrilo 2023)}\label{sec:silver-stepsizes}

The papers \cite{hedging1,hedging2} developed the \emph{silver stepsizes} in order to prove acceleration for arbitrarily large horizons $n \to \infty$. The silver stepsizes accelerate
the iteration complexity of GD from the textbook rate $O(\eps^{-1})$ to $O(\eps^{-\log_{\rho} 2})$ in the convex setting (Theorem~\ref{thm:convex}), and analogously from $O(\kappa \log \tfrac{1}{\eps})$ to $O(\kappa^{\log_{\rho} 2} \log \tfrac{1}{\eps})$ in the strongly convex setting (Theorem~\ref{thm:strongly-convex}). Here $\rho = 1+\sqrt{2}$ denotes the silver ratio. 
In these papers we conjectured that our asymptotic rate is optimal among (non-adaptive) stepsize schedules.\footnote{If one believes the optimality of our asymptotic rates, then one can ask even more fine-grained questions about the hidden constant in the big-O notation. We further conjecture that the silver stepsize schedule is exactly optimal in the strongly convex setting for the performance metric $\|x_n - x^*\|/\|x_0 - x^*\|$ that it was originally designed for. In the convex setting, an interesting line of follow-up work has improved the hidden-constant by almost a factor of $3$ and presented conjecturally optimal stepsizes for a variety of performance metrics, see the discussion in \S\ref{sec:outlook}.}

\begin{figure}
\centering
\includegraphics[width=0.36\linewidth]{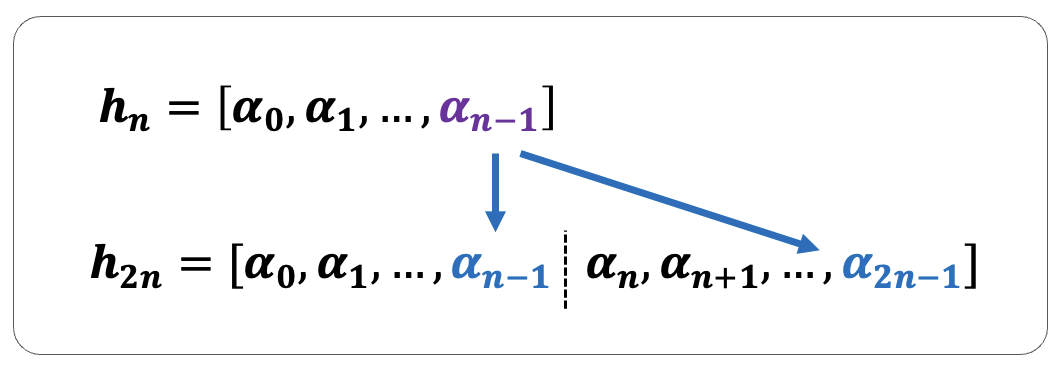}
\caption{ 
    Schematic of how the silver stepsize schedule is recursively constructed in the strongly convex setting~\cite{hedging1}. Let $h_{n}$ denote the schedule of length $n$. Then $h_{2n}$ is built from two copies of $h_{n}$, but with a key modification of the final step $\alpha_{n-1}$ in $h_n$ (purple). It is ``split'' into a shorter step $\alpha_{n-1}$ and longer step $\alpha_{2n-1}$ (blue), given by the two roots of a quadratic equation, exactly analogous to the $2$-step construction in \S\ref{sec:2-step}.}
\label{fig:silver-recursion}
\end{figure}

The silver stepsizes are constructed in a recursive manner from the $2$-step solution in Theorem~\ref{thm:2-step}, see Figure~\ref{fig:silver-recursion}. This leads to highly unconventional features: the silver stepsize schedule is non-monotonically time-varying, fractal-like, and has peaks that are exponentially infrequent but exponentially large. In the strongly-convex setting, the silver stepsizes are also approximately periodic with period of size $\kappa^{\log_{\rho} 2}$. See~\cite{hedging1} for details.

\begin{figure}
    \centering
    \includegraphics[width=0.36\linewidth]{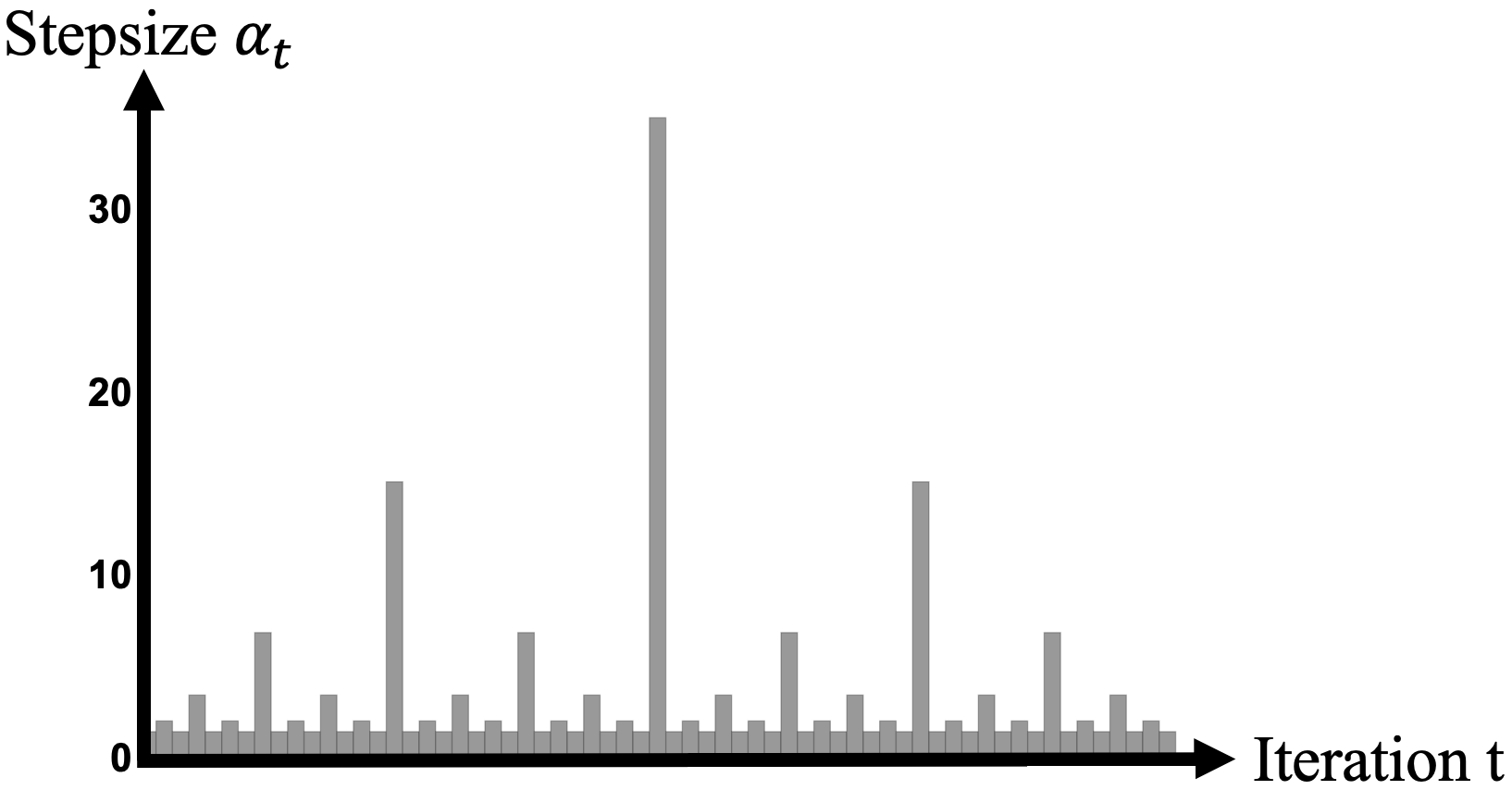}
    \caption{ 
        The silver stepsize schedule we proposed for convex optimization~\cite{hedging1,hedging2}. The first $63$ stepsizes are shown. They are highly unconventional: time-varying, non-monotonic, fractal-like, and sometimes use extremely aggressive stepsize ``spikes''
        of size $\alpha_{2^k-1} = 1+\rho^{k-1}$ where $\rho = 1+\sqrt{2}$ denotes the silver ratio.
    }
    \label{fig:silverstepsize}
\end{figure}

The silver stepsize schedule simplifies in the (non-strongly) convex setting. See Figure~\ref{fig:silverstepsize}. This is the formal limit as the strong convexity parameter $m \to 0$, or equivalently $\kappa \to \infty$. The resulting schedule has simple direct definitions, both recursively and explicitly. For simplicity, below we normalize the smoothness $M=1$ (without loss of generality since running GD on $f$ amounts to rescaling the stepsizes by $1/M$). 
\begin{itemize}
    \item \textbf{Recursive definition.} For any integer $n = 2^{k}-1$, we recursively construct the schedule $h_{2n+1}$ of length $2n+1$ from the schedule $h_{n}$ of length $n$ via
    \begin{align}
    	h_{2n+1} = [h_n, \; 1 + \rho^{k-1}, \; h_n]
    	\label{eq:steps-recursive}
    \end{align}
    with base case $h_1 = [\sqrt{2}]$. This results in the pattern $[\sqrt{2}, \; 2, \; \sqrt{2}, \; 1+\sqrt{2}, \dots]$ depicted in Figure~\ref{fig:silverstepsize}. 
    \item \textbf{Explicit definition.} The $t$-th stepsize is
    \begin{align}
        \alpha_t = 1 + \rho^{\nu(t+1)-1}
        \label{eq:steps-explicit}
    \end{align}
    where $\nu(t)$ denotes the $2$-adic valuation of $t$, i.e., the largest non-negative integer $i$ such that $2^i$ divides $t$. For example, $\nu(1) = 0$, $\nu(2) = 1$, $\nu(3) = 0$, $\nu(4) = 2$, etc. 
\end{itemize}

\paragraph*{Algorithm analysis: multi-step descent.}
The core technical challenge is that hedging requires deviating from standard analyses, which bound the progress of each step individually (e.g., using the descent lemma) and then sum the progress. Such \emph{one-step} analyses are too shortsighted to prove acceleration. Instead, one must directly analyze the  \emph{multi-step progress} in order to show that different iterations of GD synergize with each other. These holistic analyses are necessary for proving \textit{any} benefit of time-varying stepsizes. 

\par Of course, directly analyzing the long-term progress is much more difficult. This involves answering two interrelated questions, both of which are non-trivial:
\begin{itemize}
    \item \emph{Analysis question:} how to analyze a given stepsize schedule?
    \item \emph{Design question:} how to design the optimal stepsize schedule?
\end{itemize}
Our starting point is the performance estimation problem (PEP) framework, pioneered by~\cite{drori2014performance} and refined by~\cite{pesto}, which enables numerically solving the analysis question using semidefinite programming (SDP).
PEP addresses the analysis question by linearly combining valid inequalities for the class of functions under consideration; in the convex case, these are the cocoercivities relating  function values and gradients at different points.
Although PEP is helpful, it is important to emphasize that it is \textit{not} a complete solution to the problem at hand. One challenge is that PEP does not fully solve the analysis question: it only produces a numerical estimate of the convergence rate for a fixed number of iterations $n$, whereas establishing acceleration requires rigorous symbolic proofs for arbitrarily large $n$. The biggest limitation, however, is that PEP does not directly tackle the design question: finding the stepsize schedule with fastest convergence is non-convex in all existing PEP formulations\footnote{It is an interesting question if there are alternative formulations or parameterizations of the stepsize design question that are convex.}. This poses a key difficulty for discovering hedging strategies. 

Our work~\cite{hedging1,hedging2} proposes a technique called \emph{recursive gluing} which helps make both the design and analysis questions tractable for arbitrarily large $n$. In terms of the design question, as described above, we recursively construct the silver stepsize schedule of size $2n$ by gluing together two copies of the schedule of length $n$, modulo modification to a constant number of stepsizes. A key insight is that we can leverage this recursive structure of the stepsize schedule in order to also recursively approach the analysis question: we show that one can ``recursively glue'' two copies of the convergence proof for the smaller schedules, modulo modification to a constant number of inequalities. The key point is that the proof complexity does not increase in $n$. Indeed, in this way, proving the $2n$-step convergence rate from the $n$-step rate is no harder than proving the $2$-step convergence rate from the textbook $1$-step rate. See~\cite{hedging1,hedging2} for full details.

\section{Outlook}\label{sec:outlook}

This general principle of multi-step descent opens up many directions for the design and analysis of algorithms, as it suggests a potential missed opportunity for any optimization algorithm analyzed with traditional one-step analyses. In just the past three years since~\cite{hedging1,hedging2}, many papers have already followed up. For example, recent work has showcased the generality of this stepsize-hedging phenomenon by extending the results to anytime convergence~\cite{kornowski2024open,zhang2024anytime}, constrained and proximal settings~\cite{bok2024accelerating,bok2025optimized,wang2024relaxed}, Riemannian settings~\cite{park2025acceleration}, robustness and inexact gradient settings~\cite{vernimmen2024tight,vernimmen2025empirical,bai2025generalization}, operator-splitting settings~\cite{abbasza2025convergence}, min-max settings~\cite{shugart2025min,shugart2025negative}, random stepsizes~\cite{AltPar24random}, and other performance metrics~\cite{GrimmerShuWang24,zhang2024accelerated,grimmer2024composing}. The line of work~\cite{GrimmerShuWang24,wang2024relaxed,zhang2024accelerated,grimmer2024composing} has also improved the hidden constant in the big-O for the non-strongly convex setting (although still with the same asymptotic rate $\log_{\rho} 2$ as the silver stepsize schedule), and in particular~\cite{zhang2024accelerated,grimmer2024composing} developed composition/concatenation analysis frameworks which 
simplify and refine the recursive gluing technique of~\cite{hedging1,hedging2} and suggest stepsize schedules that are conjecturally minimax-optimal for several performance metrics. 
Recent work has also used large stepsizes to obtain faster rates or sharper dynamical understanding for logistic/classification losses, including separable logistic regression, regularized logistic regression, non-separable settings, and related two-layer-network models~\cite{axiotis2023gradient,wu2024large,zhang2025minimax,wu2025large,meng2024gradient,cai2024large,meng2025gradient,crawshaw2026tight,crawshaw2025constant}. The machine learning community has also been excited by intriguing similarities with similar stepsize schedules previously explored in empirical ML, e.g., the cyclical stepsize schedules of \cite{smith2019super,smith2017cyclical} as well as learned schedules in parametric optimization~\cite{SambharyaStellato,sambharya2026learning,sambharya2025learning}. Previously, time-varying stepsizes had no provable benefit even in (non-quadratic) convex optimization; this line of work provides a first step towards explaining the unreasonable effectiveness of time-varying stepsizes in empirical deep learning. But much more is needed to realize the potential of time-varying hedging as a powerful algorithmic primitive in both theory and practice.

\bibliographystyle{informs2014}
\bibliography{hedging}

\end{document}